
\input xy
\xyoption{all}

\input amssym.def
\input amssym.tex

\magnification=\magstephalf

\hfuzz=15pt

\def\im{{\frak m}}
\def\reg{\hbox{\rm reg}}

\def\indeg{\hbox{\rm indeg}}

\def\D{{\Delta}}
\def\d{{\delta}}
\def\ol{\overline}

\def\lra{\longrightarrow}
\def\fini{{$\quad\quad\square$}}

\def\fitt{\hbox{\rm Fitt}}
\def\reg{\hbox{\rm reg}}

\def\coker{\hbox{\rm coker}}

\def\indeg{\hbox{\rm indeg}}

\def\sym{\hbox{\rm Sym}}
\def\tor{\hbox{\rm Tor}}

\def\hom{\hbox{\rm Hom}}

\def\s{{\sigma}}

\font\nrm=cmcsc10 at10pt

 \font\eightit=cmti10 scaled 900
\font\bigrm=cmb10 scaled 1600


\bigskip\bigskip\bigskip\bigskip

\vskip 1cm\centerline{\bigrm Bounds for  the Castelnuovo-Mumford regularity
of  modules}

\bigskip\bigskip\bigskip

 \centerline {\bf Marc CHARDIN, Amadou Lamine
FALL, Uwe NAGEL}

\bigskip\bigskip

\centerline {\bf Abstract} \bigskip{\eightit }
  {\it We establish bounds for the
Castelnuovo-Mumford regularity of a finitely generated graded module
 and its symmetric powers in terms of the degrees of the generators of the module and the degrees of their relations. We extend to modules (and improve) the known bounds
for homogenous ideal in a polynomial ring established by Galligo, Guisti, Caviglia and Sbarra.
 }

\bigskip

\bigskip\bigskip
\centerline{\bf 1. Introduction} \medskip
Bayer and Stillman proved in [BS] that the complexity of an ideal (or a module) is the same as the one of its generic initial ideal. This connection motivated the search for bounds on the Castelnuovo-Mumford regularity in terms of the degrees of the generators of an ideal in order to control the complexity of Gr\"obner basis computations.\par
For an ideal $I$ in a polynomial ring $R=k[X_{1},\ldots,X_{n}]$, generated in degree at most $d$, Galligo ([Ga1], [Ga2]) and Giusti ([Giu]) proved the following bound if the characteristic of the field $k$ is zero:
$$ \reg(I)\leq (2d)^{2^{n-2}}.
$$
A weaker bound was given in any characteristic by Bayer and Mumford ([BM]) following an argument of Mumford ([Mum]): $$\reg(I)\leq (2d)^{(n-1)!}.$$
In [CS] Caviglia and Sbarra proved that the bound of Galligo and Giusti is valid in any characteristic. The example of Mayr and Meyer  ([MM]) shows that the doubly exponential behavior in $n$ cannot be avoided. \medskip
We extend the sharpest known  bounds  for ideals to finitely generated graded modules. As $\reg(I)=\reg(R/I)+1$, the case of ideal corresponds to cyclic modules generated in degree zero and extending the result to arbitrary modules one should take into account the degrees of the generators and of their relations.
Results on the regularity of modules have also been proved independently by Brodmann and G\"otsch in [BG]. The estimates
they provide are compared to ours in Remark 3.7. Let us only point out that the exponent in our bound depends on the dimension
of the support of the module which might be much smaller then the dimension of the ring appearing in their estimate.

We proceed in two steps. First  we establish bounds in the case of modules supported in dimension at most 1.
Second we extend  the method of Caviglia and Sbarra to modules. It allows us to proceed by induction on the dimension of the module.

For the first step we use an argument first introduced by Gruson, Lazarsfeld and Peskine in proving regularity bounds for reduced curves: a complex which is not too far from being acyclic is enough  to estimate the regularity. Our results are rather general as they hold for modules over each standard graded $R_0$-algebra where $R_0$ is an artinian local ring, provided the dimension of the modules is at most 1, and for every module over a graded Cohen-Macaulay $R_0$-algebra in higher dimension. We establish the following bounds: \medskip

{\bf Theorem.} {\it Let $R$  be a standard graded Cohen-Macaulay algebra over the artinian local ring $R_0$ and let  $M \neq 0$ be  a graded $R$-module of
finite type. Assume $M$ is generated by $n$ elements of non-negative degrees. Let $c$ and $\d$ be  the codimension and the dimension of the support of $M$ (so that $c+\d =\dim R$), respectively.  If $M$ is generated in degrees at most $B-1$ and related in degrees at most $B$, then$:$\smallskip

{\rm (i)} If $\d \leq 1$ and $c>0$,
$
\reg (M)\leq \reg (R)+(\dim R+n-1)B-\dim R
$.
\smallskip

{\rm (i)'} If $\d \leq 1$ and $c=0$,
$
\reg (M)\leq \reg (R)+B-1.
$
\smallskip

{\rm (ii)} If $\d\geq 2$ and $c>0$,
$$
\reg(M) \leq
\left[ \deg (R)(\reg (R)+(c+n)B-c){{c+n-1}\choose{c}}B^{c}\right] ^{2^{\d -2}} .
$$
\smallskip

{\rm (ii)'} If $\d\geq 2$ and $c=0$,}
$
\reg(M) \leq
\left[ n\deg (R)(\reg (R)+B)\right] ^{2^{\d -2}} .
$

\bigskip

This theorem can be used to improve the bound given by Caviglia and Sbarra in [CS] (see Remark 3.6).

It would be interesting to extend the bound in arbitrary dimension to the larger class of standard graded algebras over an artinian local ring. This could be achieved by extending  Proposition 3.3  to this more general situation. \bigskip\bigskip

 {\bf 2. Regularity bounds for modules of dimension at most 1}\bigskip
 Let $R = \oplus_{i \geq 0} R_i$ be a standard graded algebra over the artinian local ring $(R_0, {\frak n})$, i.e.
 $R_{1}$ is a finite $R_0$-module and $R = R_0 [R_1]$, and let $M$ be a graded $R$-module with a finite presentation
 $ \displaystyle F\buildrel{\varphi}\over\lra G\lra M \lra 0,$ where
 $F=\displaystyle\bigoplus_{i=1}^{m}R[-b_{i}]$, $G=\displaystyle\bigoplus_{i=1}^{n}R[-a_{i}]$ and $b_{1}\geq b_{2}\geq \ldots \geq b_{m}$.  \par

Set $\hbox{---}^{*}:=\hom_{R}(\hbox{---},R)$ and $\sigma := \sum_{i=1}^{n}a_{i}$. Let $E^{(l)}_{\bullet}$ be the complex of $R$-graded free $R$-modules associated to  $\varphi$ (see [E1]):
 $$0 \lra N^{(l)}_{m-n}[\sigma ]\buildrel{\delta} \over \lra N^{(l)}_{m-n-1}[\sigma ]\lra \ldots \lra N^{(l)}_{l}[\sigma ]\buildrel{\varepsilon} \over \lra L^{(l)}_{l}\buildrel{\nu} \over
 \lra L^{(l)}_{l-1}\lra\ldots\lra L^{(l)}_{0}\lra 0$$ with $N^{(l)}_{s}=\sym_{s-l}G^{*}\bigotimes \bigwedge^{n+s}F$  for $l\leq s\leq m-n$ and
 $L^{(l)}_{s}=\sym_{l-s}G\bigotimes\bigwedge^{s}F$ for $0\leq
 s\leq l.$
 \medskip

The complex $E^{(l)}_{\bullet}$ has graded degree zero differentials, its homology is supported in the support of $M$,
  $H_{0}(E^{(l)}_{\bullet})=\sym_{R}^{l}(M)$ for $l>0$, and $H_{0}(E^{(0)}_{\bullet})={\rm Fitt}_{0}^{R}(M)$. The complexes $E^{(0)}_{\bullet}$ and
$E^{(1)}_{\bullet}$ are  known as
Eagon-Northcott and Buchsbaum-Rim complexes, respectively. \medskip

We denote the local cohomology modules of $M$ with support in the irrelevant maximal ideal $\im = \oplus_{i > 0} R_i$ by $H^{i}_{\im}(M)$. Notice that $H^{i}_{\im}(M) \cong H^{i}_{\im + {\frak n}}(M)$.  We set
$$
a_{i}(M):=\sup\{ \mu\ \vert \ H^{i}_{\im}(M)_{\mu}\not= 0\}.
$$
Thus, $a_{i}(M):=-\infty$ if $H^{i}_{\im}(M) = 0$. Recall that $\reg (M)=\max \{ a_{i}(M)+i\}$.\medskip

 We will  provide a bound for the regularity of symmetric
  powers of $M$ in terms of the degrees of the generators of $M$ and of
  its first module of syzygies. We start by treating the case of modules
  of dimension at most 1. For simplicity, we state separately the result in case the ring has  dimension
  at most one.

\bigskip
{\bf Proposition 2.0.} {\it Let $R$  be a standard graded algebra  of dimension at most 1 over the artinian local ring $R_0$, $M$ a graded $R$-module with a finite presentation $\bigoplus_{i=1}^{m}R[-b_{i}] \buildrel{\varphi} \over \lra \bigoplus_{i=1}^{n}
R[-a_{i}] \lra M \lra 0.$ Assume that $ b_{1}\geq b_{2}\geq \cdots \geq b_{m}$, set $b:=\max\{ b_{i}\}=b_{1}$ and $a:=\max\{ a_{i}\}$. Then, for $l > 0$, \smallskip

(i) If $\dim R=0$, $\reg (R/\fitt^{0}_{R}(M))\leq \reg (R)$ and $
\reg(\sym^{l}_{R}(M))\leq \reg (R)+la$,\smallskip

(ii) If $\dim R=1$,
$$
\eqalign{
\reg(\sym^{l}_{R}(M))&\leq \max\{ a_{0}(R)+la,a_{1}(R)+(l-1)a+b\}\cr
&\leq \reg (R)+\max\{ la,(l-1)a+b-1\} ,\cr
}
$$
and
$$
\eqalign{
\reg(R/\fitt^{0}_{R}(M))
&\leq \reg (R)+\max\{ 0,\sum_{i=1}^{n}(b_{i}-a_{i})-1\} .\cr
}
$$
}

Over higher-dimensional rings we have:
\bigskip

 {\bf Theorem 2.1.} {\it Let $R$  be a standard graded algebra of dimension $d\geq 2$ over the artinian local ring $R_0$, $M \neq 0$ a graded $R$-module of dimension at most 1 with a finite presentation $M=\coker (\bigoplus_{i=1}^{m}R[-b_{i}] \buildrel{\varphi} \over \lra \bigoplus_{i=1}^{n}R[-a_{i}])$. (Note that we must have $m\geq n+d-2$.)

 Assume that $ b_{1}\geq b_{2}\geq \cdots \geq b_{m}$ and $ a_{1}\geq a_{2}\geq \cdots \geq a_{n}$. Set $D_{l}:=\sum_{i=1}^{l}(b_{i}-1)$ and $\D:=\sum_{i=1}^{\min\{ m,n+d-1\}}b_{i}-\sum_{i=1}^{n}a_{i}-(d-1)a_{n}-d$.\smallskip
 Let $M':= \coker (\bigoplus_{i=1}^{n+d-2}R[-b_{i}] \buildrel{\varphi '} \over \lra \bigoplus_{i=1}^{n}R[-a_{i}])$, where $\varphi '$ is the restriction of $\varphi$. \smallskip

Then,

(i) $
\reg(R/\fitt^{0}_{R}(M))\leq \reg (R)+\D,\
$
\smallskip

(ii) for $l \leq d-1$,
$
\reg(\sym^{l}_{R}(M))  \leq \left \{ \eqalign{ & \reg (R)+\max\{ D_{l},\D+la_{n}\}
 \; \; \; if \; M\not= M' \;  or \; l\leq d-2 \cr
& \reg (R)+D_{d-1}  \; \; \; if \; M=M' \;  and \; l = d-1,
} \right.
$
\smallskip

(iii) for $l\geq d$, $
\reg(\sym^{l}_{R}(M))\leq \reg (R)+ D_{d}+(l-d)a_{1}.
$
}
\smallskip


\bigskip We will need the following lemma to prove  the above results. It builds on ideas in [GLP] and generalizes Lemma 5.9 in [E2]. \bigskip

 {\bf Lemma 2.2.} {\it Let $C$ be a complex of finite graded $R$-modules with $C_{i}=0$ for $i<0$. If  $\dim(H_{i}(C))\leq i$ for $i>0$, then
 $$
a_{i}(H_{0}(C))\leq \max_{j\geq 0}\{a_{i+j}(C_{j})\} , \quad \forall i.
$$
In particular
$$
\reg(H_{0}(C))\leq \max_{0\leq j\leq \dim R}\{\reg (C_{j})-j\} .
$$
}\bigskip {\it Proof of  Lemma 2.2.} We consider the graded  double
complex ${\cal C}^{\bullet}_{\im}C$ where $\im$ is the maximal
homogeneous ideal of $R$ and ${\cal C}^{\bullet}_{\im}E$ is the \v Cech
complex on $E$. It gives rise to two spectral sequences. One of them
has as second terms ${'E}_{2}^{pq}=H^{p}_{\im}(H_{q}(C))$. Since
$\dim(H_{i}(C))\leq i$ for $i>0$, ${'E}_{2}^{pq}=0$ for $p>q>0$.
This implies that  ${'E}_{2}^{p0}\simeq {'E}_{\infty}^{p0}$ for each
$p$.

The other spectral sequence has as first terms, ${''E}_{1}^{pq}=H^{p}_{\im}(C_{q})$.
It follows that  $({'E}_{\infty}^{i0})_{\mu}\simeq
({'E}_{2}^{i0})_{\mu}=H^{i}_{\im}(H_{0}(C))_{\mu}$ vanishes if
$H^{p}_{\im}(C_{q})_{\mu}=0$ for $p=q+i$.\fini
\bigskip

\noindent
{\it Proof of Proposition 2.0}:\smallskip
If $\dim (R)=0$, then $\reg (R/\fitt^{0}_{R}(M))\leq  \reg (R)$ and, as $\sym^{l}_{R}(M)$ is quotient of $\sym^{l}_{R}(G)$,
$$\reg (\sym^{l}_{R}(M))\leq \reg (\sym^{l}_{R}(G))=\reg(R)+la_{1}.$$

If $\dim R=1$, then by applying Lemma 2.2 to $E^{(l)}_{\bullet}$ we get, for $l=0$ and $m \geq n$, $\reg (R/\fitt^{0}_{R}(M))\leq  \max\{ \reg (L^{(0)}_{0}),\reg ( L^{(0)}_{1})-1\}$
and, for $l>0$,
$$
\reg (\sym^{l}_{R}(M))\leq \max\{ a_0 (L^{(l)}_{0}), a_1 ( L^{(l)}_{1})-1\} \leq \reg (R)+\max\{ la, b +(l-1)a-1\}.
$$
The result follows.
\fini\medskip

\noindent
{\it Proof of Theorem 2.1}:\smallskip

Modifying the generators and relations of  $M$, the modules $F:=\bigoplus_{i=1}^{m}R[-b_{i}] $ and $G:=\bigoplus_{i=1}^{n}R[-a_{i}]$ can be decomposed into $F=F_{1}\oplus F_{2}\oplus F_{3}$ and $G=G_{1}\oplus G_{2}$ to obtain a presentation of $M$ of the following type :
$$
\xymatrix{
F_{1}\oplus F_{2}\oplus F_{3}\ar^{ \pmatrix{\psi&0\cr 0&1\cr 0&0\cr } }[rr]& &G_{1}\oplus G_{2}}
$$
where $\psi$ is a minimal presentation of $M$.

Notice that by passing to this minimal presentation, $\D$  and $D_{l}$ can only decrease.
It follows that we may assume that $ \displaystyle F\buildrel{\varphi}\over\lra G\lra M \lra 0$ is a minimal presentation
of $M$.  This assumption implies that $b_{j}>a_{n}$ for any
$j$. Notice also that $b_{j}>a_{1}$ for $j=1,\ldots ,d-1$, because $M'' :=\coker( \displaystyle F\buildrel{\varphi_1}\over\lra R[-a_{1}])$ (where $\varphi_1$ is given by the first column of $\varphi$) is a quotient of $M$, hence its dimension is at most $1$, which implies that at least $d-1$ of the degrees $b_{i}-a_{1}$ of
the entries of $\varphi_1$ are positive.

The modules $H_{s}(E^{l}_{\bullet})$ are supported in the support of $M$ (see for instance [E2, A2.59]). As
$\dim M\leq 1$ it implies that  $\dim(H_{s}(E^{l}_{\bullet}))\leq 1$ for all $s$. Therefore, by Lemma 2.2,
$$
\reg(H_{0}(E^{l}_{\bullet}))\leq \max_{0\leq s\leq d}\{ \reg (E^{l}_{s})-s\} . \leqno{(*)}
$$

The module $L^{(l)}_{s}=\sym_{l-s}G\bigotimes\bigwedge^{s}F$ is a graded free  $R$-module generated by
elements of degrees $(a_{i_{1}}+\cdots +a_{i_{l-s}})+(b_{j_{1}}+\cdots b_{j_{s}})$ with
$i_{1}\leq \cdots \leq i_{l-s}$ and $j_{1}< \cdots < j_{s}$. Hence
$$\reg (L^{(l)}_{s})=\reg (R)+\sum_{i=1}^{s}b_{i}+(l-s)a_{1} .$$

The module $N^{(l)}_{s}=\sym_{s-l}G^{*}\bigotimes\bigwedge^{n+s}F$ is a graded free  $R$-module generated by
elements of degrees $-(a_{i_{1}}+\cdots +a_{i_{s-l}})+(b_{j_{1}}+\cdots b_{j_{n+s}})$ with
$i_{1}\leq \cdots \leq i_{s-l}$ and $j_{1}< \cdots < j_{n+s}$. Hence
$$\reg (N^{(l)}_{s}[\sigma])=\reg (R)+\sum_{i=1}^{n+s}b_{i}+(l-s)a_{n}-\sigma .$$

Notice that, for $s \leq \min\{ l,d-1\},\;$
$$ \reg (L^{(l)}_{s}) \leq \reg (R)+D_{l}+s, $$
 because $b_{j}\geq a_{1}+1$ for $j\leq d-1$, and that
$$
\reg (L^{(l)}_{d}) \leq \reg (R)+ D_{d}+(l-d)a_{1} +d.
$$

In the case $M'\not= M$, one has $m\geq n+d-1$ and,  for $l \leq s \leq d-1$,
$$
\eqalign{
\reg (N^{(l)}_{s}[\sigma])&=\reg (R)+\sum_{i=1}^{n+s}b_{i}+(l-s)a_{n}-\sigma\cr
&= \reg (R)+\D+la_{n}-\sum_{i=n+s+1}^{n+d-1}b_{i}+(d-1-s)a_{n}+d\cr
&\leq \reg (R)+\D+la_{n}-(d-1-s)(b_{n+d-1}-a_{n})+d\cr
&\leq \reg (R)+\D+la_{n}+s+1\cr
}
$$
(because $b_{j}>a_{n}$ for any $j$).

A similar computation shows that $\reg (N^{(l)}_{s}[\sigma])\leq \reg (R)+\D+la_{n}+s+1$ for $s\leq m-n
=d-2$ in the case $M=M'$ (recall that in this case $E_{p}^{(l)}=0$ for $p\geq d$).

Inequality $(*)$ and the above estimates for $\reg (L^{(l)}_{s})$ and $\reg (N^{(l)}_{s}[\sigma])$ prove the  inequalities
stated in the Theorem.\fini


\bigskip

 {\bf Corollary 2.4.}
{\it Let $R$ be a standard graded algebra over a field  and let $M$ be a graded $R$-module of
dimension at most $1$. Assume $M$ is generated by $n$ elements of  degrees between $0$ and $B-1$ and related in degrees at most $B$. If $\dim R > 0$ or $n > 1$, then
 $$\reg(M)\leq \reg (R)+(\dim R+n-1)B-\dim R.$$
}

{\it Proof.}
This  immediately follows from Proposition 2.0 and Theorem 2.1 as $0\leq a_{i}\leq B-1$ and $b_{j}\leq B$ for each $i$ and $j$.
\fini\bigskip


{\bf 3.  Regularity bounds for modules of arbitrary dimension}\bigskip

 Let $R$ be a standard graded algebra over the artinian local ring $R_0$ and let  $M \neq 0$ be a graded $R$-module of dimension $\d$ presented by
 $F=\bigoplus_{i=1}^{m}R[-b_{i}] \buildrel{\varphi}\over \lra G=
\bigoplus_{i=1}^{n}R[-a_{i}] \lra M \lra 0.$ \medskip

We may write $R=S/J$, where $J$ is a graded $S$-ideal and $S$ is a polynomial ring over $R_{0}$. Assuming that $J$
has no element of degree one, this presentation is unique  and $S=\sym_{R_{0}}(R_{1})$.\medskip

Set
$$
b_{i}^{R}(M):=\sup\{ \mu \ \vert\ \tor_{i}^{R}(M,R_{0})_{\mu}\not= 0\}.
$$
Thus,   $b_{i}^{R}(M)=-\infty$ if $\tor_{i}^{R}(M,k) = 0$. Recall that $\reg (M)=\max_{i}\{ b_{i}^{S}(M)-i\}$. The initial degree of $M$ is denoted by $\indeg M = \min \{\mu \ \vert \ M_{\mu} \neq 0\}$. Furthermore, we write $\lambda (M)$ for the length of $M$ as $R_0$-module if it is finite, and we set $h^{i}_{\im}(M)_{\mu} = \lambda  (H^{i}_{\im}(M)_{\mu})$.

We will extend the technique of Caviglia and Sbarra to modules. We begin with a
result that generalizes [CS, 2.2] :\medskip

{\bf Lemma 3.1.}  {\it Let $M$ be a finitely generated graded $R$-module and let
$l$ be a linear form such that  $K :=0:_{M}(l)$ has finite length.
Set $M':=M/H^{0}_{\im}(M)$, $\ol{M}:=M/lM$, $\ol{M'}:=M'/lM'$, and $a:=a_{0}(M)-\indeg (H^{0}_{\im}(M))+1$.

Then, for all integers $\mu$, \smallskip

{\rm (i)} $
\lambda (K_{\geq \mu})=h^{0}_{\im}(M)_{\mu}+h^{0}_{\im}(\ol{M})_{>\mu}-h^{0}_{\im}(\ol{M'})_{>\mu}
\geq h^{0}_{\im}(M)_{\mu}$, \smallskip

{\rm (ii)}
$
h^{0}_{\im}(M)_{\mu +a}\leq \sum_{j=1}^{a}h^{0}_{\im}(\ol{M})_{\mu +j}-\sum_{j=1}^{a}h^{0}_{\im}(\ol{M'})_{\mu +j}
$, \smallskip

{\rm (iii)}
$
\reg(M)\leq \mu - 1 +h_{\im}^{0}(M)_{\mu}
$,
provided
$ \mu\geq \max \{ b_{0}^{R}(M) + h - 1 ,\ b_1^R (M) - 1,\ \reg (\ol{M}) + 1\}$
where $h:=\max \{ b_{0}^{S}(J),1\}$.
}\bigskip

{\it Proof.} Consider the commutative diagram with exact rows and columns
$$
\xymatrix{
 &0\ar[d]&0\ar[d]&&\\
0\ar[r]&K\ar[r]\ar[d]&K\ar[r]\ar[d]&0\ar[d]&\\
0\ar[r]&H^{0}_{\im}(M)\ar[r]\ar^{\times l}[d]&M\ar[r]\ar^{\times l}[d]&M'\ar[r]\ar^{\times l}[d]&0\\
0\ar[r]&H^{0}_{\im}(M)[1]\ar[r]\ar[d]&M[1]\ar[r]\ar[d]&M'[1]\ar[r]\ar[d]&0\\
0\ar[r]&H^{0}_{\im}(M)/lH^{0}_{\im}(M)[1]\ar[r]\ar[d]&\ol{M}[1]\ar[r]\ar[d]&\ol{M'}[1]\ar[r]\ar^{ }[d]&0\\
&0&0&0.&\\
}
$$
It shows that the length of  $K_{\mu}$ is
$$
\lambda (K_{\mu}) = h^{0}_{\im}(M)_{\mu}-h^{0}_{\im}(M)_{\mu +1}+h^{0}_{\im}(\ol{M})_{\mu +1}-h^{0}_{\im}(\ol{M'})_{\mu +1}
$$
and the first equality follows. Set $F^{j}:=0:_{M}(l^{j})$ and notice that $F^{0}=0$ and $F^{1}=K$. Using $l F^j = F^{j-1} \cap l M$ if $j \geq 2$, we get the exact sequence
$$
0\lra F^{j}/F^{j-1}\buildrel{\times l}\over{\lra}F^{j-1}/F^{j-2}[1]\lra (F^{j-1}+lM)/(F^{j-2}+lM)[1]\lra 0.
$$
It provides that
$$
\lambda (F^{j}_{\mu}/F^{j-1}_{\mu})=\lambda (F^{j-1}_{\mu +1}/F^{j-2}_{\mu +1})-\lambda ((F^{j-1}+lM)_{\mu +1}/(F^{j-2}+lM)_{\mu +1}),
$$

thus in particular, $\lambda (F^{j}_{\mu}/F^{j-1}_{\mu})\leq \lambda (K_{\mu +j-1})$ for all $j \geq 1$. As $F_{j}=F_{j+1} = H^{0}_{\im}(M)$ if $j\geq a$, it follows that
$$
\eqalign{
h^{0}_{\im}(M)_{\mu}&= \sum_{j=1}^{a}\lambda (F^{j}_{\mu}/F^{j-1}_{\mu})\cr
&\leq \sum_{j=1}^{a}\lambda (K_{\mu +j-1})\cr
&\leq h^{0}_{\im}(M)_{\mu}-h^{0}_{\im}(M)_{\mu +a}+\sum_{j=1}^{a}h^{0}_{\im}(\ol{M})_{\mu +j}-\sum_{j=1}^{a}h^{0}_{\im}(\ol{M'})_{\mu +j},\cr}
$$
which proves the second inequality.

For (iii), notice that $\reg (M) = \min \{a_0 (M), \reg (\ol{M})\}$. Suppose we know for some  $\mu \geq \reg (\ol{M})$ that $h^{0}_{\im}(M)_{j} >  h^{0}_{\im}(M)_{j+1}$ whenever $j \geq \mu$ and $h^{0}_{\im}(M)_{j} \neq 0$. Then it follows that  $a_0 (M) \leq \mu - 1 + h^{0}_{\im}(M)_{\mu}$, thus $\reg (M) \leq \mu - 1 + h^{0}_{\im}(M)_{\mu}$.

If $\mu\geq \reg (\ol{M})$ we have the exact sequence $0\lra K_{\mu}\lra H^{0}_{\im}(M)_{\mu}\lra H^{0}_{\im}(M)_{\mu +1}\lra 0$.
It implies that it suffices to estimate the degrees of the generators of $K$ as $S$-module.
Now the exact sequence $0\lra K\lra M\lra M[1]\lra \ol{M}[1]\lra 0$ provides
$$
b_0^S (K) \leq \max \{ b_{0}^{S}(M),\ b_{1}^{S}(M)-1,\ b_{2}^{S}(\ol{M})-1\} .
$$
We get a free resolution of $M$ as $S$-module from the double complex obtained by taking a free resolution of $M$ as $R$-module and then resolving each of the occurring free $R$-modules over $S$. This implies in particular that $b_1^S (M) \leq \max \{b_1^R (M),\ b_0^R (M)  + b_0^S (J) \}$. It follows that
$$
b_0^S (K) \leq \max \{ b_{0}^{R}(M) + h - 1 ,\ b_1^R (M)-1,\ \reg (\ol{M}) + 1\}.
$$
Now the last statement follows.
\fini
\medskip

Lemma 3.1 allows us to establish the following recursion:
\medskip

{\bf Lemma 3.2.} {\it Let $l_1,\ldots,l_{s+1} \in R$ be linear forms and set $M_{i}:=M/(l_{1},\ldots ,l_{i})M$, $i = 0,\ldots,s+1$.
Assume that, for each $i$,\ $K_{i}:=\ker (M_{i}\buildrel{\times l_{i+1}}\over{\lra}M_{i}[1])$ has finite length and that $M$ is generated in non-negative degrees.

Set $Q_{i}=\max\{\reg(M_{i}),\lambda(K_{i}),b_{1}^{R}(M)-2,b_{0}^{R}(M)+\max \{ 1,b_{0}^{S}(J)\} -2\}+1$ for $0\leq i\leq s$.

Then, for each $i = 0,\ldots,s-1$,\ $
Q_{i}\leq Q_{i+1}^{2}$. In particular,
$$
\reg (M)\leq Q_{s}^{2^{s}}.
$$ }

{\it Proof.} Lemma 3.1 (i) provides $h^{0}_{\im}(M_{i})_{\mu}\leq \lambda ((K_{i})_{\geq \mu}) \leq \lambda (K_i)$. Since $H^{0}_{\im}(M_{i}) \subset M_i$ does not have elements of negative degree, we get
$$
h^{0}_{\im}(M_{i}) := \lambda (H^{0}_{\im}(M_{i})) \leq (a_{0}(M_{i})+1)\cdot \lambda (K_{i}).
$$
Define $r_{i}:=\max \{ b_{1}^{R}(M)-2,\ b_{0}^{R}(M)+\max \{ 1, b_{0}^{S}(J)\} -2,\ \reg(M_{i})\}$. Then $Q_i = 1 + \max \{r_i, \lambda (K_i)\}$. Furthermore, set $R_{(i)} = R/(l_1,\ldots,l_i)R$. Using $b_0^{R_{(i)}} (M_i) \leq b_0^R (M)$ and $b_1^{R_{(i)}} (M_i) \leq b_1^R (M)$ and applying Lemma 3.1 (iii) to the $R_{(i)}$-module $M_i$, we obtain:
$$
\eqalign{
\reg(M_{i})&\leq    r_{i+1}+h^{0}_{\im}(M_{i})_{r_{i+1}+1}\cr
& \leq r_{i+1}+\lambda ((K_{i})_{> r_{i+1}})\cr
&\leq r_{i+1}+\lambda (K_{i}).}
$$
In particular, this implies
$$
r_{i}\leq r_{i+1}+\lambda (K_{i}).
$$
Using the first inequality obtained above, we conclude that
$$
\lambda (K_{i})
\leq h^{0}_{\im}(M_{i+1})\leq (a_{0}(M_{i+1})+1) \cdot \lambda (K_{i+1})
\leq (r_{i+1}+1) \cdot  \lambda (K_{i+1}). \leqno{(**)}
$$
It follows that $\lambda (K_{i})\leq Q_{i+1}(Q_{i+1}-1)$ and $r_{i}\leq (Q_{i+1}-1)+Q_{i+1}(Q_{i+1}-1)=Q_{i+1}^{2}-1$, thus $Q_i \leq Q_{i+1}^2$, as claimed.\fini
\medskip
To apply the previous lemma, we need to estimate the degree (or multiplicity) of a graded module
in terms of the degrees appearing in a graded presentation of it.\bigskip

{\bf Proposition 3.3.}  {\it Let $R$ be a standard graded Cohen-Macaulay ring and let $M$ be a graded $R$-module of codimension $c > 0$ that is generated by $n$ elements of degrees $a_{1}, \ldots ,a_{n}$ and whose  first syzygy module is generated in degrees $b_{1}\geq \ldots \geq b_{s}$. Then,  $s\geq c+n-1$ and
$$\deg(M)\leq \deg (R)\sum_{1\leq i_{1}\leq\cdots\leq i_{c}\leq n}\prod_{\ell =1}^{c}(b_{i_{\ell}+\ell -1}-a_{i_{\ell}}).$$
 }
\medskip

{\it Proof.}
Replacing the generators $g_{1},\ldots ,g_{s}$ of the first syzygy module by
$g'_{i}:=g_{i}+\sum_{j>i}h_{ij}g_{j}$, where $h_{ij}$ is a  sufficiently general polynomial
of degree $b_{i}-b_{j}$, we may assume that the $(c+n-1)\times n$ matrix $H$ corresponding
to the relations $g_{1},\ldots ,g_{c+n-1}$ has its ideal of maximal minors of codimension
$c$. Since $M$ is a quotient of the module $P$, this implies $\deg (M) \leq \deg (P)$. We will now compute the degree of $P$.

In case $R$ is a polynomial ring over a field $k$,  [Fu, 14.4.1] shows that the degree of $P$ is equal (up to the sign) to
the coefficient of order $c$ in the expansion
of
$$
{{\prod_{i=1}^{s}(1-b_{i}t)}\over{\prod_{j=1}^{n}(1-a_{j}t)}}.
$$
Now, setting $\s_{p}$ for the $p$-th symmetric function in $s$ variables, and $m_{q}$ for the
sum of monomials of degree $q$ in $n$ variables, one has:
$$
{\prod_{i=1}^{s}(1-b_{i}t)}=\sum_{p=0}^{s}(-1)^{p}\s_{p}(b_{1},\ldots ,b_{s})t^{p}
$$
and
$$
{{1}\over{\prod_{j=1}^{n}(1-a_{j}t)}}=\sum_{q\geq 0}m_{q}(a_{1},\ldots ,a_{n})t^{q}.
$$
It follows that $\deg (P)=\sum_{p+q=c}(-1)^{p-c}\s_{p}(b_{1},\ldots ,b_{s})m_{q}(a_{1},\ldots ,a_{n})$,
which can be rewritten as $$\sum_{1\leq i_{1}\leq\cdots\leq i_{c}\leq n}\prod_{\ell =1}^{c}(b_{i_{\ell}+\ell -1}-a_{i_{\ell}}).$$

Set $a:=(a_{1}, \ldots ,a_{n})$, $b:=(b_{1}, \ldots ,b_{c+n-1})$ and let $H_{a,b}(\mu )$ be the Euler-Poincar\'e characteristic of the Buchsbaum-Rim complex $(E^{(1)}_{\bullet})_{\mu}$. Since $P$ is resolved by $E^{(1)}_{\bullet}$, there is polynomial $P_{a,b}(t)$ such that
$$
\sum_{\mu}H_{a,b}(\mu )t^{\mu}={{P_{a,b}(t)}\over{(1-t)^{p-c}}}=(1-t)^{c}P_{a,b}(t){{1}\over{(1-t)^{p}}} = (1-t)^{c}P_{a,b}(t){{P_{R}(t)}\over{(1-t)^{p}}},
$$
where $p:=\dim (R)$, and $\sum_{1\leq i_{1}\leq\cdots\leq i_{c}\leq n}\prod_{\ell =1}^{c}(b_{i_{\ell}+\ell -1}-a_{i_{\ell}}) = \deg (P) = P_{a,b}(1)$. The same computations give the analogous result if $R$ is a polynomial ring over an artinian local ring $R_0$.

If $R$ is not a polynomial ring over $R_0$, then the Buchsbaum-Rim complex is still a resolution of $M$. Since  the shifts occurring in the modules $E^{(1)}_{i}$ are given by the same expressions in terms of $a$ and $b$ as above,
the Euler-Poincar\'e characteristics $H'_{a,b}(\mu )$ of $(E^{(1)}_{\bullet})_{\mu}$ of degree $\mu$ is given by
$$
\sum_{\mu}H'_{a,b}(\mu )t^{\mu}=(1-t)^{c}P_{a,b}(t){{P_{R}(t)}\over{(1-t)^{p}}}
$$
where $p:=\dim R$ and $P_{R}(1)=\deg (R)$. The conclusion follows.
\fini

\medskip

{\bf Corollary 3.4.}  {\it With the hypotheses of Proposition 3.3, set $a:=\min_{i}\{ a_{i}\}$. Then
$$
\deg (M)\leq \deg (R){{r+n-1}\choose{n-1}}\prod_{i=1}^{r}(b_{i}-a).
$$
}

{\bf Theorem  3.5.} {\it Let $R$  be a standard graded Cohen-Macaulay ring over the artinian local ring $R_0$ and let $M \neq 0$ be graded $R$-module of
finite type. Assume $M$ is generated by $n$ elements of non negative degrees. Let $c$ and $\d$ be  the codimension and the dimension of the support of $M$ (so that $c+\d =\dim R$), respectively.  If $M$ is generated in degrees at most $B-1$ and related in degrees at most $B$, then$:$\smallskip

{\rm (i)} If $\d \leq 1$ and $c>0$,
$
\reg (M)\leq \reg (R)+(\dim R+n-1)B-\dim R
$.
\smallskip

{\rm (i)'} If $\d \leq 1$ and $c=0$,
$
\reg (M)\leq \reg (R)+B-1.
$
\smallskip

{\rm (ii)} If $\d\geq 2$ and $c>0$,
$$
\reg(M) \leq
\left[ \deg (R)(\reg (R)+(c+n)B-c){{c+n-1}\choose{c}}B^{c}\right] ^{2^{\d -2}} .
$$
\smallskip

{\rm (ii)'} If $\d\geq 2$ and $c=0$,}
$$
\reg(M) \leq
\left[ n\deg (R)(\reg (R)+B)\right] ^{2^{\d -2}} .
$$

{\it Proof.} Parts (i) and (i)' are proved in Theorem 2.1 and Proposition 2.0.

Let $D:=\deg (R)$ and $r:=\reg (R)$. We may assume that the  field $R_0/{\frak n}$ is infinite.

If $\dim M=1$, one has $\lambda (K)=\lambda (M/lM)-\deg (M)$ where $l \in R$ is any linear form such that $K:=0:_{M}(l)$ has finite length. Corollary 3.4 applied to $M/lM$ provides
$\lambda (M/lM)\leq D{{c+n-1}\choose{n-1}}B^{c}$ if $c>0$. Notice that $\lambda (M/lM)\leq nD$ if $c=0$. It follows that $\lambda (K)<D{{c+n-1}\choose{n-1}}B^{c}$ if $\dim M=1$ and $c>0$ and $\lambda (M/lM)\leq nD$ if $\dim M=1$ and $c=0$.

Let $\d =2$ and $c>0$ and choose  $l$ as a general element in $R_1$. Notice that $b_{0}^{R}(M)+\max\{ 1,b_{0}^{S}(J)\}\leq (B-1)+(r+1)$. Hence, using the notation of Lemma 3.2, part (i) implies
$$
r_1 = \max\{ b_{1}^{R}(M) -2, b_{0}^{R}(M)+\max\{ 1,b_{0}^{S}(J)\} -2, \reg (M/l M)\} \leq r+(c+n)B-(c+1).
$$
Thus we get from Lemma 3.1 (i) and (iii) that $\reg (M)< \lambda (K)+r+(c+n)B-c$.
Moreover, inequality $(**)$ in the proof of Lemma 3.2 provides
$$
\lambda (K) \leq (r+(c+n)B-c)\left[ D{{c+n-1}\choose{n-1}}B^{c}-1\right].
$$
It follows that
$$
\eqalign{
\max \{ \reg (M),\lambda (K)\}+1&\leq (r+(c+n)B-c)\left[ D{{c+n-1}\choose{n-1}}B^{c}-1\right] +r+(c+n)B-c\cr
&= D(r+(c+n)B-c){{c+n-1}\choose{n-1}}B^{c}.
\cr
}
$$

If $\d =2$ and $c=0$, the estimates are respectively $r_1\leq r+B-1$ by part (i)', $\reg (M)< \lambda (K)+r+B$ by Lemma 3.1 (i) and (iii) and $\lambda (K) \leq (r+B)(nD-1)$ by  inequality $(**)$ in the proof of Lemma 3.2.
It follows that $\max \{ \reg (M),\lambda (K)\}+1\leq nD(r+B)$ in this case.\medskip

Finally, if $\d >2$, the result follows by induction using Lemma 3.2 with $s=\d -2$.\fini\medskip

{\bf Example 3.6.} Consider the special case where $R = S$ is a polynomial ring over an artinian local ring $R_0$ and  $M=S/I$ is of codimension $c$, dimension $\d$, and $I$ is generated in degree at most $B$. Then Theorem 3.5 provides
$$
\reg(I) \leq
\left[ (c+1)B^{c+1}\right] ^{2^{\d -2}}.
$$
Notice that if $\dim S=p \geq 2$ and $c=1$, then $I=(F)J$, where $F$ has degree $e\leq B$ and the ideal $J$ is
generated in degree at most $B-e$ and has codimension $c' \geq 2$. Then we get
$$
\reg(I) =e+\reg (J)\leq
e+\left[ (c'+1) (B-e)^{c'+1}\right] ^{2^{p-4}}.
$$
It follows that the regularity of every ideal in $I \subset S:=R_0 [X_{1},\ldots ,X_{p}]$ that is generated in degree at most $B$ is
bounded above by $p(B-1)+1$ if $p\leq 3$, and by $\left[ 3B^{3}\right] ^{2^{p-4}}$ if $p\geq 4$. Notice that these results also
follow from [Ch1, 3.3] or [Sj, Theorem 2] and [CF, Theorem 1] by using Lemma 3.2. In fact, it gives the slightly refined bound (which also follows from the proof of Theorem 3.5):
$$
\reg (I) \leq \left[ 3B^{2}(B-1)\right] ^{2^{p-4}}+1, \quad {\rm provided} \ p \geq 4.
$$
This in turn improves the bound of
Caviglia and Sbarra [CS]: $\reg (I) \leq \left[ B^{2}+2B-1\right] ^{2^{p-3}}$.\medskip

{\bf Remark 3.7.}(i) In part (iii) of Theorem 3.5 the same arguments show a refined inequality if $c>0$. Namely, if
$M$ is generated in degrees $a_{1},\ldots ,a_{s}$ and related in degrees $b_{1}\geq \cdots \geq b_{s}$, then
(unless $s=c+n-1$ in which case $M$ is Cohen-Macaulay of regularity $\reg (R)+\sum_{j}b_{j}-\sum_{i}a_{i}-c\min\{ a_{i}\}$)
one has
$$
\reg(M) \leq
\left[ \deg (R)(\reg (R)+b_{1}+\cdots +b_{c+n}-c)\sum_{1\leq i_{1}\leq\cdots\leq i_{r}\leq n}\prod_{\ell =1}^{r}(b_{i_{\ell}+\ell -1}-a_{i_{\ell}})\right] ^{2^{\d -2}} .
$$
\bigskip

(ii) In [BG, 6.3], Brodmann and G\"otsch prove the following bound:
$$
\reg (M)\leq \left[ \reg (R)+(n+1)\deg (R)+B+1\right] ^{2^{p-1}},
$$
where $p = \dim R$. In fact, they give a bound that is a little more precise (using more data on the $a_{i}$'s and $b_{j}$'s). The main difference to our bound is that in our estimate the exponent depends on
the dimension of the module as opposed to the dimension of the ambient ring. Furthermore, our bound  is slightly stronger, even in
the case of a module supported in small codimension. \medskip

{\bf Remark 3.8.} If $M$ is generated in degrees between $0$ and $a$ and related in degrees at most $B$, then $\sym^{l}_{R}(M)$ is generated in degrees between $0$ and $la$ and related in degrees at most $B+(l-1)a$. Applying
Theorem 3.5 in this situation, one obtains if $a\leq B-1$ and $\d =\dim (M)\geq 2$,
$$
\reg(\sym^{l}_{R}(M)) \leq
\left[ \deg (R)(\reg (R)+(c+n)(B+(l-1)a-c){{c+n-1}\choose{c}}(B+(l-1)a)^{c}\right] ^{2^{\d -2}} .
$$

{\bf Remark 3.9.}
 It would be interesting to extend the bound in Theorem 3.5 to the class of arbitrary standard graded algebras over an artinian local ring.

\bigskip\bigskip
\centerline{\bf REFERENCES }
\bigskip

\noindent [BM] D. Bayer, D. Mumford, {\it What can be computed in algebraic geometry?,} Computational algebraic geometry and commutative algebra, Symposia Mathematica, {\bf XXXIV} (1993), 1-48 \medskip
\noindent [BS] D. Bayer, M. Stillman, {\it A criterion for detecting $m$-regularity,} Invent. Math. {\bf 87} (1987), 1-11.\medskip
\noindent [BG] M. Brodmann, T. G\"otsch, {\it Bounds for the Castelnuovo-Mumford regularity}, Preprint 25--2005, Institut f\"ur Mathematik, Universit\"at Z\"urich.\medskip
\noindent [CS] G. Caviglia, E. Sbarra, {\it Characteristic-free bounds
for the Castelnuovo-Mumford regularity,} Compositio Math. {\bf 141} (2005), 1365-1373. \medskip
 \noindent [Ch1] M. Chardin,  {\it Regularity of
ideals and their powers},  Pr{\'e}publication {\bf 364},  Institut de
math{\'e}\-matiques de Jussieu, Mars 2004.\medskip
 \noindent [Ch2] M. Chardin,  {\it Some results and questions on Castelnuovo-Mumford regularity},
Pr{\'e}publication.\medskip
\noindent [CF] M. Chardin, A. L. Fall, {\it Sur la r\'egularit\'e de Castelnuovo-Mumford des id\'eaux, en dimension deux},
C. R. Acad. Sci. Paris {\bf 341} (2005), 233-238. \medskip
\noindent [E1] D. Eisenbud, { \it
Commutative Algebra with a View Toward Algebraic Geometry}, Graduate
Texts in Mathematics {\bf 150}, Springer-Verlag, New York, 1995.\medskip
\noindent [E2] D. Eisenbud, { \it The Geometry of Syzygies},
Graduate Texts in Mathematics {\bf 229}, Springer-Verlag, New York,
2005.\medskip
\noindent [Fu] W. Fulton, { \it Intersection Theory},
Second Edition,  Springer-Verlag, New York,
1998.\medskip
\noindent [Ga 1] A. Galligo, {\it A propos du th\'eor\`eme de pr\'eparation de Weierstrass. Fonctions de plusieurs variables complexes}, Lectures Notes in Mathematics  {\bf  409}   (1974), 543-579.\medskip

\noindent [Ga 2] A. Galligo, {\it  Th\'eor\`eme de division et stabilit\'e en g\'eometrie analytique locale,} Ann. Inst. Fourier {\bf 29}  (1979), 107-184.\medskip

\noindent [Giu] M. Giusti, {\it Some effectivity problems in polynomial ideal theory,} Eurosam {\bf 84}, Lect. Notes Comp. Sci {\bf 1984}, 159-171. \medskip

\noindent [GLP] L. Gruson, R. Lazarsfeld, and C. Peskine, {\it On a theorem of Castelnuovo, and the equation defining space curves,} Invent. Math. {\bf 72} (1983), 491-506. \medskip

\noindent [MM] E. Mayr, A. Meyer, {\it The complexity of the word problem for commutative semigroups and polynomial ideals}, Adv. Math. {\bf 46} (1982), 305-329. \medskip

\noindent [Mum] D. Mumford, {\it Lectures on algebraic curves on algebraic surfaces}, Princeton University Press, Princeton, New Jersey, 1966. \medskip

\noindent [Sj] R. Sj\"ogren, {\it On the regularity of graded $k$-algebras of Krull dimension $\leq 1$},  Math. Scand.  {\bf 71}  (1992),   167--172.\medskip

\bigskip\bigskip\bigskip

{\nrm Marc Chardin}\smallskip

Institut de Math\'ematiques de Jussieu,

CNRS \& Universit\'e Pierre et Marie Curie,

4, place Jussieu, 75005 Paris,

 France\smallskip

{\it e-mail:} chardin@math.jussieu.fr

\bigskip

{\nrm Amadou Lamine Fall}\smallskip

Facult{\'e} des sciences, Universit{\'e} Cheikh Anta Diop

Dakar

S{\'e}n{\'e}gal\smallskip

{\it e-mail:} fall@math.jussieu.fr

\bigskip

{\nrm Uwe Nagel}\smallskip

Department of Mathematics,

University of Kentucky,

715 Patterson Office Tower,

Lexington, Kentucky 40506-0027,

USA \smallskip

{\it e-mail:} uwenagel@ms.uky.edu

\end